\documentclass[11pt]{amsart}
\usepackage{graphicx}
\usepackage{color}
\usepackage{latexsym}
\usepackage{amscd}
\usepackage{a4wide}
\usepackage{epsfig}
\usepackage{float}
\usepackage{amsmath, amssymb, amsthm, amsfonts, amsgen} 
\usepackage[outdir=./]{epstopdf}

\graphicspath{ {./} }

\usepackage{esint}

\usepackage{enumitem}
\usepackage{cite}

\usepackage{color}

\renewcommand{\[}{\begin{equation*}}
\renewcommand{\]}{\end{equation*}}
\newcommand{\grad}{\ensuremath{\nabla}}
\newcommand{\p}{\ensuremath{\partial}}
\newcommand{\ve}{\ensuremath{\varepsilon}}

\newcommand{\la}{\ensuremath{\lambda}}
\renewcommand{\th}{\ensuremath{\theta}}

\renewcommand{\ss}{\ensuremath{\subset}}

\newcommand{\ra}{\ensuremath{\rightarrow}}

\newcommand{\mb}[1]{\ensuremath{\;\mbox{#1}\;}}
\newcommand{\cd}{\ensuremath{\cdot}}
\newcommand{\ti}{\ensuremath{\times}}

\newcommand{\cof}{\textnormal{cof}\;}
\renewcommand{\div}{\textnormal{div}\;}

\newcommand{\R}{\mathbb{R}}

\newcommand{\N}{\mathbb{N}}

\newcommand{\A}{\mathcal{A}}
\newcommand{\F}{\mathcal{F}}

\newcommand{\Om}{\Omega}

\newtheorem{de}{Definition}[section]
\newtheorem{thm}[de]{Theorem}

\newtheorem{re}[de]{Remark}
\allowdisplaybreaks
\begin{document}
\title[(High frequency-) uniqueness criteria in in- and compressible elasticity]{(High frequency-) uniqueness criteria for $p-$growth functionals in in- and compressible elasticity} 
\author[Marcel Dengler]{M. Dengler}
\address[M. Dengler]{Fliederweg 1, 72189 Vöhringen, Germany.}
\email{marci.dengler@web.de}


\maketitle

\begin{abstract}
\noindent In this work our main objective is to establish various (high frequency-) uniqueness criteria. Initially, we consider $p-$Dirichlet type functionals on a suitable class of measure preserving maps $u: B\ss \R^2 \to \R^2,$ $B$ being the unit disk, and subject to suitable boundary conditions.
In the second part we focus on a very similar situations only exchanging the previous functionals by a suitable class of $p-$growing polyconvex functionals and allowing the maps to be arbitrary.In both cases a particular emphasis is laid on high pressure situations, where only uniqueness for a subclass, containing solely of variations with high enough Fourier-modes, can be obtained.
\end{abstract}


\section{Incompressible setting and results}
\label{sec:4.2.5}
Let $B\ss \R^2$ be the unit ball. For $u_0\in L^p(B,\R^2)$ with $2\le p<\infty$ and $\det\grad u_0=1$ a.e. in $B$ we define
\[
\A_{u_0}^{p,c}:=\{u\in W^{1,p}(B,\R^2): \det\grad u=1 \;\mb{a.e. in}\;B ,\;u-u_0\in C_c^\infty(B,\R^2)\}
\]
and for every $u\in \A_{u_0}^{p,c}$ we define
\begin{equation}E(u)=\int\limits_B{f(x,\grad u)\;dx}\label{eq:HFU.I.2.1},\end{equation}
where $f$ is given by
\[f(x,\xi)=\nu(x)|\xi|^p \]
for a.e. $x\in B$ and $\xi\in \R^{2\ti2}.$ 
Moreover,  $\nu\in L^\infty(B)$ is supposed to satisfy $\nu(x)\ge0$ a.e.\! in $B.$ Since $\nu$ is allowed to take on the value $0,$ the integrand could indeed disappear for some $x \in B.$ Additionally, $f(x,\cd)$ is convex in its 2nd variable for a.e. $x\in B.$ All of the above, is making $E$ a version of the $p-$Dirichlet functional. \vspace{0.5cm}

As usual, we are interested in the  corresponding minimization problem
\begin{equation}
\inf\limits_{u\in\A_{u_0}^{p,c}} E(u).
\label{eq:USPS:1.1}
\end{equation}
Notice, that the missing uniformity might cause troubles guaranteeing the existence of a minimizer, however, this is without any consequence for this work.
Assuming for now, that the situation is such that the minimum is indeed obtained and there are some corresponding minimizing maps or more generally corresponding stationary points of $E,$ which are defined as follows:

\begin{de}(Stationary point)
We say that $u$ is a stationary point of $E(\cdot)$ if there exists a function $\la$, which we shall henceforth refer to as a pressure, belonging to $W^{1,1}(B)$ and such that
\begin{align}\label{def:SP1} \div( \nabla_{\xi}f(x,\nabla u) + p \la(x) \, \cof \nabla u) = 0 \quad \mathrm{in} \ \mathcal{D}'(B).
\end{align}
\end{de}

 Now we can state the main result for the incompressible scenario. Recall, that for any vector $y=y_Re_R+y_\th e_\th\in\R^2$ we define its maximum norm via $|y|_{\infty}:=\max\{y_R,y_\th\}.$

\begin{thm} [High frequency uniqueness] Let $2\le p<\infty,$ assume $u_0\in L^p(B,\R^2)$ to be the boundary conditions and let $u\in \A_{u_0}^{p,c}$ be a stationary point of $E,$ as given in \eqref{def:SP1}. 
Furthermore, let $\sigma(x):=\sqrt{\nu(x)|\grad u(x)|^{p-2}}\in L^{\frac{4}{p-2}}(B)$ and assume that there exists $l\in \N$ s.t.
\begin{equation}|\sigma,_\th(x)|\le l \sigma(x) \mb{for a.e.}\;x\in B\label{eq:HFU.UC.G.01} \end{equation}
holds.\\

Then the following statements are true:

i) \textbf{(purely high modes.)}  Suppose the corresponding pressure $\la$ exists and satisfies
\begin{equation}|\grad \la(x)R|_{\infty}\le \frac{n}{\sqrt{2}}\nu(x)|\grad u|^{p-2}\mb{for a.e.}\;x\in B\label{eq:HFU.IC.G.02}\end{equation}
for some $n\in\N.$\\

Then $u$ is a minimizer of $E$ in the subclass 
\[\F_{n_*}^{p,\sigma,c}=\left\{v\in \A_{u_0}^{p,c}|\eta=v-u\in W_0^{1,p}(B,\R^2) \;\mb{and}\; \sigma\eta=\sum\limits_{j\ge n}(\sigma\eta)^{(j)} \right\},\]
where $n_*:=n+l.$
Moreover, if there exists a constant $\sigma_0>0$ s.t. $\sigma(x)\ge\sigma_0>0$ for any $x\in B$ and inequality \eqref{eq:HFU.IC.G.02} is strictly satisfied on a non-trivial set, then $u$ is the unique minimizer in $\F_{n_*}^{p,\sigma,c}$.\\

ii) \textbf{($0-$mode and high modes.)} Suppose the corresponding pressure $\la$ exists and satisfies
\begin{equation}|\grad \la(x)R|_{\infty}\le \frac{\sqrt{3}m\nu(x)|\grad u|^{p-2}}{2\sqrt{2}}\mb{for a.e.}\;x\in B\label{eq:HFU.IC.G.03}\end{equation}
for some $m\in\N.$\\
Then $u$ is a minimizer of $E$ in the subclass 
\[\F_{0,m_*}^{p,\sigma,c}=\left\{v\in \A_{u_0}^{p,c}|\eta=v-u\in W_0^{1,p}(B,\R^2) \;\mb{and}\; \sigma\eta=(\sigma\eta)^{(0)}+\sum\limits_{j\ge m_*}(\sigma\eta)^{(j)} \right\},\]
where $m_*=m+l.$
Moreover, if there exists a constant $\sigma_0>0$ s.t. $\sigma(x)\ge\sigma_0>0$ for any $x\in B$ and inequality \eqref{eq:HFU.IC.G.03} is is strictly satisfied on a non-trivial set, then $u$ is the unique minimizer in $\F_{0,m_*}^{p,\sigma,c}$.
\label{thm1.1}
\end{thm}\vspace{0.5cm}

\section{Compressible setting and results}
 For $u_0\in L^p(B,\R^2)$ with $2\le p\le\infty$ we define the set of admissible maps by
\[
\A_{u_0}^p=\{u\in W^{1,p}(B,\R^2): \; u-u_0\in C_c^\infty(B,\R^2)\}.
\]
Here we consider energies given by
\begin{equation}I(u)=\int\limits_B{\Phi(x,\grad u)\;dx},\label{eq:HFU.UC.G.0}\end{equation}
where the integrand is of the form
\begin{equation}\Phi(x,\xi)=\frac{\nu(x)}{p}|\xi|^p+\Psi(x,\xi,\det\xi),\label{eq:HFU.UC.GA}\end{equation}
Now we want $(\xi,d)\mapsto \Psi(x,\xi,d)$ to be convex for a.e.\! $x \in B,$ making $\Psi(x,\cd)$ a convex representative of a polyconvex function a.e. in $B$. 
Moreover, the function $\nu\in L^\infty(B),$ satisfying $\nu(x)\ge0$ a.e.\! in $B,$ is ought to be optimal, in the sense, that there can be no term of the form  $a(x)|\xi|^q$ for any $q\ge p$ in $\Psi.$ However, $\Psi$ might be negative. Finally, we want $\Phi$ to be of $p-$growth,  supposing that there exits $C\in L^\infty(B)$ with $C(x)\ge0$ a.e. in $B$ s.t.
\[0\le \Phi(x,\xi)\le \frac{C(x)}{p}(1+|\xi|^p)\;\mb{for all} \;\xi\in \R^{2\ti2}\;\mb{and a.e.} \; x\in B. \]
All of the above combined guarantees that $x\mapsto \Phi(x,\grad u(x))\in L^1(B)$ for any $u\in W^{1,p}(B,\R^2)$ and hence the corresponding energy is finite. \vspace{0.5cm}

The corresponding minimization problem
\begin{equation}
\inf\limits_{u\in\A_{u_0}^p} I(u)
\label{eq:USPS:1.1}
\end{equation}
might lack a minimizer in general. Assuming, that such a minimizer/stationary point of $I$ exists, motivates the following definition: 
\begin{de}(Stationary point)
We say that $u\in\A_{u_0}^p$ is a stationary point of $I(\cdot)$ if $u$ satisfies the Euler-Lagrange equation given by
\begin{align}\label{def:SP.2}
\div( \nu(x)|\grad u|^{p-2}\grad u+\p_\xi\Psi(x,\grad u,d_{\grad u})+\p_d\Psi(x,\grad u,d_{\grad u})\cof\grad u) = 0\; \mathrm{in} \;\mathcal{D}'(B).
\end{align}
\end{de} \vspace{-0.3cm}
Now we can give an analogous result for these types of integrands.

\begin{thm} [High frequency uniqueness] Let $B\ss \R^2$ be the unit ball. Furthermore, let $2\le p\le\infty,$ assume $u_0\in L^p(B,\R^2)$ to be the boundary conditions and let $u\in \A_{u_0}^p$ be a stationary point of $I,$ as given in \eqref{def:SP.2}. Furthermore, let $\sigma(x):=\sqrt{\nu(x)|\grad u(x)|^{p-2}}\in L^{\frac{4}{p-2}}(B)$ and assume that there exists $l\in \N$ s.t.
\begin{equation}|\sigma,_\th(x)|\le l \sigma(x) \mb{for a.e.}\;x\in B\label{eq:HFU.UC.G.01} \end{equation}
holds.\\

Then the following statements are true:

i) \textbf{(purely high modes.)}  Assume there exists $n\in\N$ s.t.
\begin{equation}|\grad_x \p_d\Psi(x,\grad u(x),d_{\grad u(x)})R|_{\infty}\le \frac{n}{\sqrt2}\nu(x)|\grad u|^{p-2}\mb{for a.e.}\;x\in B.\label{eq:HFU.UC.G.1} \end{equation}

Then $u$ is a minimizer of $E$ in the subclass 
\[\F_{n_*}^{p,\sigma}=\left\{v\in \A_{u_0}^p|\;\eta=v-u\in W_0^{1,p}(B,\R^2) \;\mb{and}\; \sigma\eta=\sum\limits_{j\ge {n_*}}(\sigma\eta)^{(j)} \right\},\]
where $n_*:=n+l.$
Moreover, if there exists a constant $\sigma_0>0$ s.t. $\sigma(x)\ge\sigma_0>0$ for any $x\in B$ and inequality \eqref{eq:HFU.UC.G.1} is strictly satisfied on a non-trivial set, then $u$ is the unique minimizer in $\F_{n_*}^{p,\sigma}$.\\

ii) \textbf{($0-$mode and high modes.)} Assume there exists $m\in\N$ s.t.
\begin{equation}|\grad_x \p_d\Psi(x,\grad u(x),d_{\grad u(x)})R|_{\infty}\le \frac{\sqrt{3}m\nu(x)|\grad u|^{p-2}}{2\sqrt{2}}\mb{for a.e.}\;x\in B.\label{eq:HFU.UC.G.2}\end{equation}
Then $u$ is a minimizer of $E$ in the subclass 
\[\F_{0,m_*}^{p,\sigma}=\left\{v\in \A_{u_0}^p|\;\eta=v-u\in W_0^{1,p}(B,\R^2) \;\mb{and}\; \sigma\eta=(\sigma\eta)^{(0)}+\sum\limits_{j\ge m_*}(\sigma\eta)^{(j)} \right\},\]
where $m_*:=m+l.$
Moreover, if there exists a constant $\sigma_0>0$ s.t. $\sigma(x)\ge\sigma_0>0$ for any $x\in B$ and inequality \eqref{eq:HFU.UC.G.2} is strictly satisfied on a non-trivial set, then $u$ is the unique minimizer in $\F_{0,m_*}^{p,\sigma}$.
\label{thm:HFU.UC.G.1}
\end{thm}

Here, we contribute to John Ball's objectives, as formulated in \cite[§ 2.6]{BallOP}, which tries to get an understanding of questions related to uniqueness in elastic situations. We recently started this discussion in \cite{D1}, where we presented a uniqueness criteria, in the incompressible case\footnote{Recall, that we call situations, where the energies remain finite either incompressible elastic, if the considered admissible maps must be measure-preserving, otherwise we call it compressible elastic. Moreover, we call a model (fully) non-linear elastic, if the considered integrand $f,$ ignoring any other dependencies, satisfies $f(\xi)=+\infty$ for any $\xi\in \R^{n\ti n}$ s.t. $\det\xi\le 0$ and $f(\xi)\ra+\infty$ if $\det\xi\ra 0^+$ or $\det\xi\ra +\infty.$ Notice, that in the latter model some of the consider energies might be infinite.}, for quadratic uniformly convex integrands $f(x,\xi)$ and subject to suitable boundary conditions. Here we go beyond the latter first by allowing the more general $p-$Dirichtlet integrands but also by discussing high-pressure situations, where uniqueness can only be guaranteed for a subclass of variations which consist of purely high-modes. We also provide the analogous results for polyconvex-type functionals, which are of $p-$growth, in the compressible case. The work, which is certainly most relevant to ours is the one by Sivaloganathan and Spector \cite{SS18}. In particular, in \cite[Theorem 4.2]{SS18} they consider the same functionals $I$, as given in (\ref{eq:HFU.UC.G.0}-\ref{eq:HFU.UC.GA}), however, on a suitable class of admissible maps satisfying the constraint $\det\grad u>0 \;\mb{a.e.}$ and subject to suitable boundary data. Assuming then that $u$ is a weak solution satisfying the corresponding equilibrium equation and the condition given by

\begin{equation}|\p_d\Psi(x,\grad u,d_{\grad u})R|\le \nu(x)|\grad u|^{p-2}\mb{for a.e.}\;x\in \Om.
\label{eq:Uni.SivSP.1}
\end{equation}
Then $u$ must be a global minimizer of $I.$ Additionally, if the latter inequality is in some sense strictly satisfied then $u$ must be the unique one.  It is crucial to realise, that the criterion given by Sivaloganathan and Spector differs from ours. Indeed, \eqref{eq:Uni.SivSP.1} involves a 1st order derivative of $\Psi$ while in \eqref{eq:HFU.UC.G.2} a 2nd order derivative is used.  \\
 
A possible incomplete list of results regarding the compressible setting might look as follows: Initially, it is well known, that uniformly convex functionals possess unique global minimizers, see for instance \cite[§ 3.3]{Kl16}. Knorps and Stuart showed in \cite{KS84}, that for a strongly quasiconvex integrand defined on a star-shaped domain and subject to linear boundary data $u_0=Ax,$ any $C^2-$stationary point needs to agree with $Ax$ everywhere. A generalisation can be found in \cite{T03}. These results have been transferred to the incompressible by Shahrokhi-Dehkordi and Taheri  \cite{ST10} and the fully non-linear case by Bevan \cite{B11}. In \cite{C14}, Cordero presented a uniqueness result guaranteeing a unique minimizer for strongly quasiconvex $C^2-$integrands if the given boundary data is smooth and small enough.
 Zhang \cite{Z91} discusses situations, in in- and in compressible ones, where the considered energies are polyconvex and the under pure displacement boundary conditions and subject to that the Jacobian must be strictly positive a.e. Then the corresponding minimizer must agree with the solution of the corresponding Euler-Lagrange equation, which is highly non-trivial due to the weak spaces involved and the lack of compactness of the constraint.
 On the contrary, non-uniqueness for minimizers of strongly polyconvex functionals has been established by Spadaro in \cite{S08}. However, these counterexamples rely highly on the fact that the determinant can take on negative values, which is neither possible in the incompressible nor in the NLE-stetting. John \cite{J72} and Spector and Spector \cite{SS19} obtain uniqueness of equilibrium solutions for small enough strains and under various boundary conditions. In sharp contrast, Post and Sivaloganathan \cite{PS97} construct multiple equilibrium solutions in finite elasticity. \\

A first treatment of uniqueness in incompressible elasticity can be found in \cite[Section 6]{KS84}. Much research in the incompressible setting is concerned with the double covering problem, that is, given the Dirichlet Energy $\mathbb{D}(\xi)= |\xi|^2/2$ on the unit ball in $2d$ and subject to double covering boundary conditions given by $u_2=(\cos(2\th),\sin(2\th)),$ first considered by Ball \cite{B77}.  Since then a lot of progress towards a solution has been made and partial results are available. For instance, Bevan \cite{JB14} showed that $u_2$ is the unique global minimizer up to the first Fourier-mode. This paper is central for several reasons, firstly, it introduces the type of uniqueness argument we provide here and in \cite{D1}. Moreover, it also contains the concept of high-frequency uniqueness. In \cite{BeDe21} Bevan and Deane obtained that $u_2$ is the unique global minimizer for either purely inner 
or purely outer variations and local minimality is shown for a subclass of variations allowing a certain mixture of both types. In contrast, in \cite{BeDe20}, equal energy stationary points of an inhomogeneous uniformly convex functional $(x,\xi)\mapsto f(x,\xi)$ depending discontinuously on $x$ are constructed. It remains unknown for now if these stationary points are actually global minimizers.\\


In the fully non-linear case Bevan and Yan \cite{BY07} show, that the famous BOP$-$map, constructed by Bauman, Owen and Phillips \cite{BOP91MS}, is the unique global minimizer in a suitable sub-class of admissible maps. \\

Finally, many papers address uniqueness questions, in situations where the reference domain agrees with an annulus, see for instance \cite{J72,PS97,T09,MT17,MT19,BK19}.  \\

\vspace{3mm}

\textbf{Plan of the paper:} After introducing the most important notation, which we will use in this paper, Then the proof of theorem \ref{thm1.1} will be given in §.\ref{sec:3}  followed by some important remarks in the incompressible situation. §.\ref{sec.4} then discusses the compressible case instead and, in particular, a proof of theorem \ref{thm:HFU.UC.G.1}.\\

\textbf{Notation:} For any $2\times2-$matrix $A$ we define its cofactor by
\begin{equation}
\cof A=\begin{pmatrix}a_{22} &-a_{21}\\-a_{12} & a_{11}\end{pmatrix}.
\label{eq:1.3}
\end{equation}
Moreover, we will make use of the shorthand $d_A:=\det A.$\\

The Fourier-representation for any $\eta\in C^\infty(B,\R^2)$ (For members of Sobolev- spaces one might approximate) is given by
\begin{align*}
\eta(x)=\sum\limits_{j\ge0}\eta^{(j)}(x),\;\mb{where}\;\eta^{(0)}(x)=\frac{1}{2}A_0(R), \;A_0(R)=\frac{1}{2\pi}\int\limits_{0}^{2\pi}{\eta(R,\th)\;d\th}\end{align*} 
and for any $ j\ge1$ we have
\begin{align*}
 \eta^{(j)}(x)=A_j(R)\cos(j\th)+B_j(R)\sin(j\th),
 \end{align*}
  where
\begin{align*}
A_j(R)=\frac{1}{2\pi}\int\limits_{0}^{2\pi}{\eta(R,\th)\cos(j\th)\;d\th}\;\mb{and}\;
B_j(R)=\frac{1}{2\pi}\int\limits_{0}^{2\pi}{\eta(R,\th)\sin(j\th)\;d\th}.
\end{align*}
Further we use $\tilde{\eta}:=\eta-\eta^{(0)}.$\\

\section{The incompressible case}
We immediately start with the proof in the incompressible setting. The proof is obtained by comparing energies and gaining a lower bound by means of the ELE and a Poincaré-type inequality.\\

\label{sec:3}
\textbf{Proof of Theorem \ref{thm1.1}:}\\

i) Let $u\in\A_{u_0}^{p,c}$ be a stationary point of $E$ and let $v\in\F_{n_*}^{p,\sigma,c}$ be arbitrary and set $\eta:=v-u\in W_0^{1,2}(B,\R^2)$ and $ \sigma\eta=\sum\limits_{j\ge {n_*}}(\sigma\eta)^{(j)}$ with $n_*=n+l$ assuming wlog. $\eta\in C_c^\infty(B,\R^2)$ and $\sigma\in C^\infty(B).$

We start by the standard expansion
\begin{align}
E(v)-E(u)=&\int\limits_B{\nu(x)(|\grad u+\grad \eta|^{p}-|\grad u|^{p})\;dx}&\nonumber\\
\ge&\frac{p}{2}\int\limits_B{\nu(x)|\grad u|^{p-2}|\grad\eta|^2+p\nu(x)|\grad u|^{p-2}\grad u\cd\grad\eta\;dx},&
\label{eq:HFU.IC.G.04}
\end{align}
where we used the following inequality\footnote{see, \cite[Prop A.1]{SS18}with $\sigma=0.$}\vspace{-5mm}

\begin{align}\frac{1}{p}|b|^{p}\ge\frac{1}{p}|a|^{p}+|a|^{p-2}a(b-a)+\frac{1}{2}|a|^{p-2}|b-a|^{2}.\label{eq:3.1}\end{align}

Henceforth, we shall denote the rightmost term in \eqref{eq:HFU.IC.G.04} by \[H(u,\eta):=\int\limits_B{p\nu(x)|\grad u|^{p-2}\grad u\cd\grad\eta\;dx}.\]
Recall, the ELE \eqref{def:SP1} is given by
\begin{align}
\int\limits_B{p\nu(x)|\grad u|^{p-2}\grad u\cd\grad\eta\;dx}=&-\int\limits_B{p\la \cof\grad u\cd \grad\eta\;dx}\;\mb{for all}\;\eta\in C_c^\infty(B,\R^2).
\label{eq:ELE.2}
\end{align}
In order for us, to control $H$ from below, we start by rewriting said term via the relation $\det \grad \eta=-\cof\grad u\cd\grad\eta$ a.e. and lemma 2.1(v) of \cite{D1} in the following way
\begin{align}
H(u,\eta)=&-\frac{p}{2}\int\limits_B{R((\cof\grad\eta)\grad\la)\cd\eta\;\frac{dx}{R}}&\label{eq:HFU.UC.G.21.a}\\
=&-\frac{p}{2}\int\limits_B(\la,_RR e_R+ \la,_\th e_\th)\cd\left[(\tilde{\eta}_1\tilde{\eta}_{2,\th}-\tilde{\eta}_2\tilde{\eta}_{1,\th})\frac{e_R}{R}\right.&\nonumber\\
&+\left. (\tilde{\eta}_2\eta_{1,R}-\tilde{\eta}_1\eta_{2,R})e_\th\right]\;\frac{dx}{R}.&\nonumber
\end{align}
An application of Hölder's inequality in $\R^2$, that is, for any $y,z\in \R^2$ it holds $|y\cd z|\le|y|_{\infty}|z|_{1},$  yields
\vspace{1mm}
\begin{small}
\begin{align}
H(u,\eta)\ge&-\frac{p}{2}\int\limits_B\!|R\grad\la(x)|_{\infty}\!\cd\!\left[|\eta_1\eta_{2,\th}-\eta_2\eta_{1,\th}|\frac{1}{R}+ |\eta_2\eta_{1,R}-\eta_1\eta_{2,R}|\right]\!\frac{dx}{R}.&
\end{align}
\end{small}
By applying \eqref{eq:HFU.UC.G.1} and Cauchy-Schwarz inequalities with weight $\ve=\frac{n}{\sqrt{2}}$ we get 
\begin{align*}
H(u,\eta)\ge-\frac{p}{2}\frac{n}{\sqrt{2}}&\Big[\frac{n}{\sqrt{2}}\|\sigma\eta_1\|_{L^2(dx/R^2)}^2+\frac{n}{\sqrt{2}}\|\sigma\eta_2\|_{L^2(dx/R^2)}^2&\\
&+\frac{1}{\sqrt{2}n}\Big(\|\sigma\eta_{1,\th}\|_{L^2(dx/R^2)}^2+\|\sigma\eta_{2,\th}\|_{L^2(dx/R^2)}&\\
&+\|\sigma\eta_{1,R}\|_{L^2(dx)}^2+\|\sigma\eta_{2,R}\|_{L^2(dx)}^2\big)\Big].&
\end{align*}
Using \eqref{eq:HFU.FE.2} and collecting terms yields
\begin{align*}
H(u,\eta)\ge&-\frac{p}{2}\left[\|\sigma\eta,_\th\|_{L^2(B,\R^2,\frac{dx}{R^2})}^2+\frac{1}{2}\|\sigma\eta,_R\|_{L^2(B,\R^2,dx)}^2\right]&\nonumber\\
\ge&-\frac{p}{2}\int\limits_B{\nu(x)|\grad u|^{p-2}|\grad\eta|^2\;dx},&\nonumber
\end{align*}
completing the proof. For the last step we made used of the following version of the Fourier-estimate given by
\begin{equation} n^2\int\limits_B{\sigma^2|\eta|^2\;\frac{dx}{R^2}} \le \int\limits_B{\sigma^2|\eta_{,\th}|^2\;dx}.\label{eq:HFU.FE.2}\end{equation}
This is indeed true, for this sake, first assume that  $\sigma\in C^\infty(B).$ 
We will make use of a Poincare type inequality, first established in \cite{JB14},
\begin{equation}\int\limits_B{R^{-2}|\xi_{,\th}|^2\;dx}\ge N^2\int\limits_B{R^{-2}|\xi|^2\;dx},\label{eq:Uni.HFU.1}\end{equation}
which holds true for any $\xi\in C^\infty(B,\R^2)$ if $\xi$ only consists of Fourier-modes $N$ or higher.\\

Then an application of \eqref{eq:Uni.HFU.1}, the product rule, Minkowski's inequality and \eqref{eq:HFU.UC.G.01} yields   
\begin{align} n_*\|\sigma\eta\|_{L^2\left(\frac{dx}{R^2}\right)}\le&\|(\sigma\eta),_\th\|_{L^2\left(\frac{dx}{R^2}\right)}&\nonumber\\
\le& \|\sigma\eta,_\th\|_{L^2\left(\frac{dx}{R^2}\right)}+ \|\sigma,_\th\eta\|_{L^2\left(\frac{dx}{R^2}\right)}&\nonumber\\
\le& \|\sigma\grad\eta\|_{L^2(dx)}+l \|\sigma\eta\|_{L^2\left(\frac{dx}{R^2}\right)}.&
\label{eq:HFU.FE.2a}
\end{align}
Absorbing the rightmost term of the latter expression into the LHS yields \eqref{eq:HFU.FE.2}. Additionally, \eqref{eq:HFU.FE.2a} justifies its own upgrade by remaining valid for any $\sigma\in L^2(B),$ satisfying \eqref{eq:HFU.UC.G.01}.\\

ii) Let $u\in\A_{u_0}^{p,c}$ be a stationary point of $E$ and let $v\in\F_{0,m_*}^{p,\sigma,c}$ be arbitrary and set $\eta:=v-u\in W_0^{1,2}(B,\R^2)$ and $ \sigma\eta=(\sigma\eta)^{(0)}+\sum\limits_{j\ge {m_*}}(\sigma\eta)^{(j)}$ with $m_*=m+l$ assuming wlog. $\eta\in C_c^\infty(B,\R^2)$ and $\sigma\in C^\infty(B).$

Notice, that
\begin{align}
E(v)-E(u)
\ge&\frac{p}{2}\int\limits_B{\nu(x)|\grad u|^{p-2}|\grad\eta|^2\;dx}+H(u,\eta),&
\label{eq:HFU.IC.G.05}
\end{align}
remains valid with the same mixed term
\begin{equation}H(u,\eta)=\int\limits_B{p\nu(x)|\grad u|^{p-2}\grad u\cd\grad\eta\;dx}.\label{eq:Uni.SPC.1}\end{equation}
From here on we have to argue more along the lines of the proof of \cite[thm 1.2]{D1}. Again, by the ELE \eqref{def:SP1}, the identity $\det \grad \eta=-\cof\grad u\cd\grad\eta$ a.e. and by \cite[lem 2.1.(vi)]{D1} we get
\[H(u,\eta)=-p\int\limits_B{(\cof\grad\eta^{(0)}\grad\la(x))\cd\tilde{\eta}\;dx}-p\int\limits_B{(\cof\grad\eta\grad\la(x))\cd\tilde{\eta}\;dx}=:(I)+(II).\]
Now by noting that the $0-$mode is only a function of $R,$ we get
\begin{align*}
(\cof\grad\eta^{(0)}\grad\la(x))\cd\tilde{\eta}=\frac{\la,_\th}{R}(\eta_{1,R}^{(0)}\tilde{\eta}_2-\eta_{2,R}^{(0)}\tilde{\eta}_1).
\end{align*}
Instead of just $\la,_\th$ on the right hand side of the latter equation we would like to have the full gradient of $\la.$ This can be achieved by using the basic relations $e_\th\cd e_\th=1$ and  $e_R\cd e_\th=0$ to obtain 
\begin{align*}
(\cof\grad\eta^{(0)}\grad\la(x))\cd\tilde{\eta}=(\la,_RR e_R+\la,_\th e_\th)\cd(\eta_{1,R}^{(0)}\tilde{\eta}_2-\eta_{2,R}^{(0)}\tilde{\eta}_1)\frac{e_\th}{R}.
\end{align*}
Arguing similarly for (II), and a short computation shows
\begin{align}
H(u,\eta)=&-p\int\limits_B(\la,_RR e_R+\la,_\th e_\th)\cd\left[(\tilde{\eta}_1\tilde{\eta}_{2,\th}-\tilde{\eta}_2\tilde{\eta}_{1,\th})\frac{e_R}{R}\right.&\nonumber\\
&+\left. (\tilde{\eta}_2(\eta_{1,R}^{(0)}+\eta_{1,R})-\tilde{\eta}_1(\eta_{2,R}^{(0)}+\eta_{2,R}))e_\th\right]\;\frac{dx}{R}.&
\label{eq:Uni.SPC.2}\end{align}
By Hölder's inequality in $\R^2$ we get\begin{align*} 
H(u,\eta)\ge&-p\int\limits_B|\grad \la(x)R|_\infty \left[\left|\tilde{\eta}_1\tilde{\eta}_{2,\th}-\tilde{\eta}_2\tilde{\eta}_{1,\th}\right|\frac{1}{R}\right.&\\
&+\left.\left|\tilde{\eta}_2(\eta_{1,R}^{(0)}+\eta_{1,R})-\tilde{\eta}_1(\eta_{2,R}^{(0)}+\eta_{2,R})\right|\right]\;\frac{dx}{R}.&
\end{align*}
By $|\grad \la(x)R|_\infty\le\frac{\sqrt{3}m\sigma^2(x)}{2\sqrt{2}}$ and a weighted Cauchy-Schwarz inequality, we see
\begin{align*}
H(u,\eta)\!\ge&-\frac{\sqrt{3}mp}{4\sqrt{2}}\!\left[2a\|\sigma\tilde{\eta}_1\|_{L^2\left(\frac{dx}{R^2}\right)}^2\!+2a\|\sigma\tilde{\eta}_2\|_{L^2\left(\frac{dx}{R^2}\right)}^2\!+\frac{1}{a}\|\sigma\eta_{2,R}^{(0)}+\sigma\eta_{2,R}\|_{L^2(dx)}^2\right.&\\
&\left.+\frac{1}{a}\|\sigma\tilde{\eta}_{2,\th}\|_{L^2\left(\frac{dx}{R^2}\right)}^2+\frac{1}{a}\|\sigma\eta_{1,R}^{(0)}+\sigma\eta_{1,R}\|_{L^2(dx)}^2+\frac{1}{a}\|\sigma\tilde{\eta}_{1,\th}\|_{L^2\left(\frac{dx}{R^2}\right)}^2\right].&
\end{align*}

Applying the Cauchy-Schwarz inequality, inequality \eqref{eq:HFU.FE.2}, and combining some of the norms yields
\begin{align*}
H(u,\eta)\ge-\frac{\sqrt{3}mp}{4\sqrt{2}}&\left[\left(\frac{2a}{m^2}+\frac{1}{a}\right)\|\sigma\tilde{\eta},_\th\|_{L^2\left(\frac{dx}{R^2}\right)}^2+\frac{2}{a}\|\sigma\eta,_R^{(0)}\|_{L^2(dx)}^2\right.&\\
&\left.+\frac{2}{a}\|\sigma\eta,_R\|_{L^2(dx)}^2\right].&
\end{align*}
By using $\tilde{\eta},_\th=\eta,_\th$ and $\|\sigma\eta,_R^{(0)}\|_{L^2(dx)}^2\le\|\sigma\eta,_R\|_{L^2(dx)}^2$ we obtain
\begin{align*}
H(u,\eta)\ge&-\frac{\sqrt{3}mp}{4\sqrt{2}}\left[\left(\frac{2a}{m^2}+\frac{1}{a}\right)\|\sigma\eta,_\th\|_{L^2\left(\frac{dx}{R^2}\right)}^2+\frac{4}{a}\|\sigma\eta,_R\|_{L^2(dx)}^2\right].&
\end{align*}
Choosing $a=\frac{\sqrt{3}m}{\sqrt{2}}$ and combining the norms even further yields
\begin{align*}
H(u,\eta)\ge-\frac{p}{2}\int\limits_B{\nu(x)|\grad u|^{p-2}|\grad\eta|^2\;dx},
\end{align*}
 completing the 2nd part of the proof.
\qed

\begin{re} 
 1. Notice, that \eqref{eq:HFU.UC.G.01} is especially satisfied if $p=2$ and $\nu(x)=\nu(R).$ So condition \eqref{eq:HFU.UC.G.01}
can be thought of as a natural extension of this fact to the case, where $p$ might be arbitrary and $\sigma(x)$ depends on $x$ instead of $R.$\\
2. Despite the fact that the sets $F_{n_*}^{p,\sigma,c}$ and $F_{0,m_*}^{p,\sigma,c}$ depend on $\sigma$ it remains true that if $n_*=0$ or $m_*\in\{0,1\}$ one gets uniqueness in the full class $\A^{p,c}.$ Indeed, there are two cases to consider. Firstly, let $n_*=0$ or $m_*\le1$ s.t. $n=0$ or $m=0.$ Then $\grad\la\equiv0$ and by  \eqref{eq:HFU.IC.G.04} and \eqref{eq:ELE.2} one obtains
\begin{align*}
E(v)-E(u)\ge\int\limits_B{\frac{\nu(x)}{2}|\grad u|^{p-2}|\grad\eta|^2\;dx},\end{align*}
implying, that $u$ is a global minimizer in the full class $\A^{p,c}.$ Assuming, additionally, that  $\sigma(x)\ge\sigma_0>0$ for  a.e. $x\in B$ then one can conclude that it needs to be the unique one. In the case when $m_*=1$ and $m=1,l=0,$ then  $\sigma_{,\th}=0$ and \eqref{eq:HFU.FE.2} can be applied with $m_*=m=1$ completing the argument. 
\label{re:1}
\end{re}\vspace{0.5cm}

\section{Compressible uniqueness criterion}
Here we outline the proof for the analogous compressible results, where the proof strategy remains the same.\\

\label{sec.4}
\textbf{Proof of Theorem \ref{thm:HFU.UC.G.1}:}\\
i) Let $u\in\A_{u_0}^p$ be a stationary point of $I$ and let $v\in\F_{n_*}^{p,\sigma}$ be arbitrary and set $\eta:=v-u\in W_0^{1,2}(B,\R^2)$ and $ \sigma\eta=\sum\limits_{j\ge {n_*}}(\sigma\eta)^{(j)}$ with $n_*=n+l$ assuming wlog. $\eta\in C_c^\infty(B,\R^2)$ and $\sigma\in C^\infty(B).$ Again, by the standard expansion, the subdifferential inequality for $\Psi$, and inequality \eqref{eq:3.1} we obtain
\begin{align}
I(v)-I(u)=&\int\limits_B{\frac{\nu(x)}{p}(|\grad u+\grad \eta|^{p}-|\grad u|^{p})}&\nonumber\\
&+{\Psi(x,\grad u+\grad \eta,\det\grad u+\grad \eta)-\Psi(x,\grad u,\det\grad u)\;dx}&\nonumber\\
\ge&\int\limits_B{\frac{\nu(x)}{2}|\grad u|^{p-2}|\grad\eta|^2+\nu(x)|\grad u|^{p-2}\grad u\cd\grad\eta}&\nonumber\\
&+{\p_\xi\Psi(x,\grad u,\det\grad u)\cd\grad\eta+\p_d\Psi(x,\grad u,\det\grad u)\cof\grad u\cd\grad\eta\;dx}.&\nonumber
\end{align}

Then by applying the ELE \eqref{def:SP.2} we get
\begin{align}
I(v)-I(u)\ge&\int\limits_B{\frac{\nu(x)}{2}|\grad u|^{p-2}|\grad\eta|^2+\p_d\Psi(x,\grad u,\det\grad u)d_{\grad\eta}\;dx}.&
\label{eq:HFU.UC.G.21}
\end{align}
Replacing $\la$ by $\p_d\Psi$ in the argument given in the 1st part of the proof of theorem \ref{thm:HFU.UC.G.1}, see §.\ref{sec:3}, completes the argument.\\

ii) If $\sigma\eta$ can be such that the first Fourier-mode $(\sigma\eta)^{(0)}\not=0$ then \eqref{eq:HFU.UC.G.21} does still hold, but one needs to conclude similarly to the 2nd part of the proof of theorem \ref{thm:HFU.UC.G.1}.\qed 
\begin{re}
1. Arguing as in Remark \ref{re:1} in §.\ref{sec:3}, assuming, additionally, that $\sigma(x)\ge\sigma_0>0$ for a.e. $x\in B$ then it remains true that if $n_*=0$ or $m_*\in\{0,1\}$ one gets uniqueness in the full class $\A^{p}.$\\
2. It seems reasonable to believe, that similar uniqueness criteria could be given in many other elastic scenarios if the considered situation is such, that the describing functional decomposes into two parts, the main one and the perturbative one acting like the pressure. If the perturbation then is assumed to be small in some sense there might be a chance of obtaining uniqueness.
\end{re}
\vspace{0.5cm}

\textbf{Data Declaration}
Data sharing not applicable to this article as no datasets were generated or analysed during the current study.\\

\textbf{Acknowledgments}

The author is in deep dept to Jonathan J. Bevan, Bin Cheng, and Ali Taheri for vital discussions and suggestions, the fantastic people with the Department of Mathematics at the University of Surrey and the Engineering \& Physical Sciences Research Council (EPRSC), which generously funded this work.\\


\begin{thebibliography}{10}

\bibitem{BallOP}
John~M. Ball.
\newblock {\em {Some Open Problems in Elasticity}}, chapter I.1, pages 3--59.
\newblock Springer-Verlag, NY, 2002.


\bibitem{B77}
John~M. Ball.
\newblock {Convexity conditions and existence theorems in nonlinear elasticity}.
\newblock {\em Arch. Rat. Mech. Anal.}, 64:337--403, 1977.


\bibitem{BOP91MS}
P.~Bauman, N.~C. Owen, and D.~Phillips.
\newblock {Maximal smoothness of solutions to certain
  Euler{\textendash}Lagrange equations from nonlinear elasticity}.
\newblock {\em Proc. R. Soc. Ed. Sec. A.}, 119(3-4):241--263, 1991.

\bibitem{B11}
Jonathan~J. Bevan.
\newblock Extending the Knops-Stuart-Taheri technique to $C^1$ weak local minimizers in nonlinear elasticity.\newblock {\em Proc. Amer. Math. Soc.}, 139:1667-1679, 2011.

\bibitem{JB14}
Jonathan~J. Bevan.
\newblock On double-covering stationary points of a constrained {D}irichlet
  energy.
\newblock {\em Annales de l{\textquotesingle}Institut Henri Poincare (C) Non
  Linear Analysis}, 31(2):391--411, 2014.

\bibitem{BeDe20}
Jonathan~J. Bevan and Jonathan H.~B. Deane.
\newblock {A continuously perturbed Dirichlet energy with area-preserving
  stationary points that `buckle' and occur in equal-energy pairs}.
\newblock {\em Non. Dif. Eq. App.}, 28(1), dec 2020.

\bibitem{BeDe21}
Jonathan~J. Bevan and Jonathan H.~B. Deane.
\newblock {Energy minimizing N-covering maps in two dimensions}.
\newblock {\em Var. Par. Dif. Eq. }, 60(4), 2021.

\bibitem{D1}
Jonathan~J. Bevan and Marcel Dengler.
\newblock A uniqueness criterion and a counterexample to regularity in an
  incompressible variational problem.
\newblock {\em arxiv:2205.07694}, 2022.

\bibitem{BK19}
Jonathan~J. Bevan and S.~Käbisch. 
\newblock{Twists and shear maps in nonlinear elasticity: explicit solutions and vanishing Jacobians}, 
 \newblock{\em Proc. R. Soc. Ed. Sec. A.}, 2019.
 
 \bibitem{BY07}
Jonathan~J. Bevan and X.~Yan.
\newblock Minimizers with topological singularities in two dimensional
  elasticity.
\newblock {\em {ESAIM}: Control, Optimisation and Calculus of Variations},
  14(1):192--209, sep 2007.
 
\bibitem{C14}
Judith~Campos Cordero.
\newblock{Regularity and Uniqueness in the Calculus of
  Variations}
 \newblock{\em  Ph.D. thesis, Oxford University}, 2014.
  
\bibitem{J72}
F.~John.
\newblock {Uniqueness of non-linear elastic equilibrium for prescribed boundary
  displacements and sufficiently small strains}.
\newblock {\em Commun. Pure Appl. Math.}, 25:617–634., 1972.

\bibitem{Kl16}
A.~Klaiber.
\newblock {Variationsrechnung}.
\newblock {\em Lecture Notes, Uni Konstanz}, 2016.


\bibitem{KS84}
R.~J. Knops and C.~A. Stuart.
\newblock {Quasiconvexity and uniqueness of equilibrium solutions in nonlinear
  elasticity}.
\newblock {\em Arch. Rat. Mech. Anal.}, 86(3):233–249, 1984.


\bibitem{MT17}
C.~Morris and A.~Taheri.
\newblock Twist maps as energy minimisers in homotopy classes: Symmetrisation
  and the coarea formula.
\newblock {\em Nonlinear Analysis: Theory, Methods {\&} Applications},
  152:250--275, 2017.

\bibitem{PS97}
K.~D.~E. Post and J.~Sivaloganathan.
\newblock On homotopy conditions and the existence of multiple equilibria in
  finite elasticity.
\newblock {\em Proc. R. Soc. Ed. Sect. A}, 127(3):595–614, 1997.

\bibitem{ST10}
S.~Shahrokhi-Dehkordi and A.~Taheri.
\newblock {Quasiconvexity and Uniqueness of Stationary Points on a Space of
  Measure Preserving Maps}.
\newblock {\em Journal of Convex Analysis}, 17:69–79, 2010.

\bibitem{S86}
Jeyabal Sivaloganathan.
\newblock {Uniqueness of Regular and Singular Equilibria for Spherically
  Symmetric Problems of Nonlinear Elasticity}.
\newblock {\em Arch. Rat. Mech. Anal.}, 96:97--136., 1986.

\bibitem{SS18}
Jeyabal Sivaloganathan and Scott~J. Spector.
\newblock {On the Uniqueness of Energy Minimizers in Finite Elasticity}.
\newblock {\em J. Elas.}, 133(1):73--103, feb 2018.

\bibitem{S08}
E.~N. Spadaro.
\newblock {Non-Uniqueness of Minimizers for Strictly Polyconvex Functionals}.
\newblock {\em Arch. Rat. Mech. Anal.}, 193(3):659--678, 2008.

\bibitem{SS19}
Daniel~E. Spector and Scott~J. Spector.
\newblock {Uniqueness of Equilibrium with
  Sufficiently Small Strains in Finite Elasticity}.
\newblock {\em Arch. Rat. Mech. Anal.}, 233:409--449, 2019.


\bibitem{T03}
Ali Taheri. 
\newblock {Quasiconvexity and uniqueness of stationary points in the
  multi-dimensional Calculus of Variations}
   Proc. Amer. Math. Soc.
  \textbf{131} (2003), 3101--3107,
\newblock {\em Proc. Amer. Math. Soc.}, 131:3101--3107, 2003.


\bibitem{T09}
Ali Taheri.
\newblock {Minimizing the Dirichlet Energy over a space of measure preserving
  maps}.
\newblock {\em Topological Methods in Nonlinear Analysis Journal of the Juliusz
  Schauder Center}, 33:170–204., 2009.

\bibitem{MT19}
Ali Taheri and Charles Morris.
\newblock {On the Uniqueness and Monotonicity of Energy Minimisers in the
  Homotopy Classes of Incompressible Mappings and Related Problems}.
\newblock {\em Journal of Mathematical Analysis and Applications},
  473(18):1--26, 2019.

\bibitem{Z91}
Kewei Zhang.
\newblock Energy minimizers in nonlinear elastostatics and the implicit
  function theorem.
\newblock {\em Archive for Rational Mechanics and Analysis}, 114(2):95--117,
  1991.

\end{thebibliography}
\end{document}